\documentclass{article}
\usepackage{amsfonts}
\usepackage{amsmath,amsthm}

\usepackage[all]{xy}
\usepackage{graphicx}

\usepackage{amssymb}
\usepackage{latexsym}

\DeclareMathAlphabet\EuFrak{U}{euf}{m}{n}	
\SetMathAlphabet\EuFrak{bold}{U}{euf}{b}{n}	

\newcommand{\ra}{\rightarrow}

\newcommand{\hra}{\hookrightarrow}

\newcommand{\wa}{\widehat}
\newcommand{\wt}{\widetilde}

\newcommand{\sC}{{\it C*}-}
\newcommand{\bC}{{\mathbb C}}

\newcommand{\bT}{{\mathbb T}}

\newcommand{\bZ}{{\mathbb Z}}

\newcommand{\bN}{{\mathbb N}}
\newcommand{\bP}{{\mathbb P}}
\newcommand{\bS}{{\mathbb S}}

\newcommand{\ud}{{{\mathbb U}(d)}}
\newcommand{\sud}{{{\mathbb {SU}}(d)}}

\newcommand{\eps}{\varepsilon}

\newcommand{\mA}{\mathcal A}
\newcommand{\mB}{\mathcal B}

\newcommand{\mE}{\mathcal E}

\newcommand{\mO}{\mathcal O}
\newcommand{\mP}{\mathcal P}

\newcommand{\mZ}{\mathcal Z}

\newcommand{\mG}{\mathcal G}

\newcommand{\zro}{C(X^\rho)}

\newcommand{\oro}{\mO_\rho}
\newcommand{\oroeps}{\mO_{\rho,\eps}}

\newcommand{\ii}{\iota,\iota}

\newcommand{\ers}{E^r,E^s}
\newcommand{\rhors}{\rho^r , \rho^s}

\newcommand{\hrs}{H^r , H^s}

\newcommand{\mQG}{\mathcal{QG}}
\newcommand{\mNG}{\mathcal{NG}}

\newcommand{\coe}{\mO_E}

\newcommand{\cog}{\mO_G}

\newtheorem{thm}{Theorem}[section]
\newtheorem{cor}[thm]{Corollary}
\newtheorem{lem}[thm]{Lemma}
\newtheorem{prop}[thm]{Proposition}

\newtheorem{defn}[thm]{Definition}

\theoremstyle{definition}
\newtheorem{ex}{Example}[section]

\theoremstyle{remark}
\newtheorem{rem}{Remark}[section]

\numberwithin{equation}{section}

\begin{document}

\author{           
 {\sc Roberto Conti}  \\ 
 {\footnotesize Mathematics}                                              \\
 {\footnotesize School of Mathematical and Physical Sciences}             \\
 {\footnotesize University of Newcastle, Callaghan, NSW 2308, Australia}  \\ [1mm]
 {\footnotesize Roberto.Conti@newcastle.edu.au}
\and
 {\sc Ezio Vasselli}\\
 {\footnotesize Dipartimento di Matematica}\\
 {\footnotesize University of Rome "La Sapienza"}  \\       
 {\footnotesize P.le Aldo Moro 2, I-00185 Roma, Italy}    \\[1mm]
 {\footnotesize vasselli@mat.uniroma2.it}}

\date{\today{}}

\title{Extension of automorphisms to \sC crossed products with non-trivial centre}
\maketitle

\begin{abstract}
Given 
a quasi-special endomorphism $\rho$ of 
a $C^*$-algebra
$\mA$ with nontrivial center,
we study an extension problem for automorphisms of $\mA$
to a minimal cross-product $\mB$ of $\mA$ by $\rho$.
Exploiting some aspects of the underlying
generalized Doplicher-Roberts duality theory 
based on Pimsner algebras, an obstruction 
to the existence of such extensions
is found and described in terms of sections of a suitable group bundle.
\end{abstract}


\section{Introduction}
Symmetries of $C^*$-algebras are provided by 
automorphisms.
Thus the structure of the group of
automorphisms of a $C^*$-algebra often reveals interesting informations 
about the algebra itself.
If two $C^*$-algebras are linked by some natural map, e.g. a 
homomorphism, 
then one can try to ``transfer'' 
automorphisms back and forth using this map.
In this paper we 
deal with a typical extension problem for automorphisms that 
loosely speaking takes the following form.
(We consider only the case of a single automorphism, 
but the generalization to the case of an action of, say, a locally compact group is straightforward.)
Suppose that $\mathfrak A$ is a 
C$^*$-subalgebra of the 
$C^*$-algebra $\mathfrak B$, 
and $\beta$ is 
an outer automorphism of $\mathfrak A$.
Then a simple question is
whether $\beta$ admits an extension to an automorphism of $\mathfrak B$. 
This question is appealing as it is both natural and useful in a number of situations.
Under more specific circumstances, it may even become important to require the extension to satisfy additional properties,
such as to commute with another given action of a locally compact group $\alpha$ on $\mathfrak B.$
In operator algebras, such kind of problems have already 
appeared in this or similar formulations
(also involving endomorphisms or anti-automorphisms) in several contexts,  
but we will not attempt to give a comprehensive list of all the related literature .
For the applications we have in mind, we only refer the reader to \cite{DR89A,BDLR,CoDA} and references therein.
As we explain below, the $C^*$-algebras we are interested in are suggested by the Doplicher-Roberts duality theory.

Motivated by structural problems in quantum field theory, 
Doplicher and Roberts have deeply investigated an abstract duality theory for compact groups
involving ``extensions of $C^*$-algebras by tensor $C^*$-categories'' 
which, in more concrete terms, take the form of cross-products by certain 
categories whose objects are 
semigroups of unital endomorphisms \cite{DR89A}.
More recently, some steps have been taken towards a generalization of their work 
in which
the relevant tensor $C^*$-categories are not required to have the property that the self-intertwiners of the monoidal unit reduce to the complex scalars $\bC$,
thus aiming at a more general duality theory for ``loop group-like'' objects \cite{Vas,Vas2,Vas3}.
A very important ingredient of the DR-analysis was provided by the Cuntz algebras, as $C^*$-algebras generated
by a Hilbert space. In a similar fashion, in the new setup a similar role is now taken by the Pimsner algebras, as $C^*$-algebras generated by a Hilbert bimodule \cite{Pim97}.

In this paper we start with a unital $C^*$-algebra with nontrivial center, equipped with an endomorphism that satisfies a weakened condition of permutation symmetry and in addition is quasi-special. Such endomorphisms naturally arise in the
framework of (generalized) DR duality. 
Due to the wide generality of the situation at hand the existence of 
a minimal cross-product of $\mA$ by $\rho$ 
is not automatically guaranteed, but it is in several particular 
situations,
as the one studied by H.~Baumg\"artel and F.~Lled\'o (\cite{BL}).  
When it is possible to construct 
such a cross-product (the ``field algebra''), 
which then clearly contains $\mA$, our extension problem applies.
The automorphisms that we wish to extend are those that are somewhat ``compatible'' with the given data.
Making use of the cross-product structure and appealing to some aspects of 
the DR duality theory
we reduce the problem to a more elementary one living in a Pimsner algebra.
We then tackle the problem of extending certain automorphisms of a fixed point subalgebra of the Pimsner algebra
under a group bundle action to the whole algebra, possibly in equivariant way,
and here the obstruction pops up. It should be pointed out that this 
phenomenon is genuinely new as the obstruction vanishes in the more traditional setting based on Cuntz algebras.

In order to keep things as simple as possible, in this paper
we have considered only the 
case of free bimodules as some salient features
of the underlying structure already show up in this situation, 
although in the future it would certainly make sense to consider more general bimodules. 
Also, we have focused our attention on the so-called equivariant automorphisms and 
not more general families of automorphisms, as a proper treatment
in the latter case
would probably require to replace 
the cross-product by a single endomorphism
with more general types of cross-products that are not yet available.
Finally, we stress that one might also try to elucidate whenever possible some additional properties of group bundles which might play a role in this context (e.g., a natural guess would be a `fibered version' of quasi-completeness,
cf. \cite{CoDA}).



We conclude with few words on the adopted terminology and conventions.
In this paper 
all the $C^*$-algebras will be unital.
Also, automorphism (resp. endomorphism, homomorphism) will stand for $*$-automorphism 
(resp. unital $*$-endomorphism, unital $*$-homomorphism).

\section{Background}
\label{intro}

Let $d \in \bN$, $d \geq 2$.
We denote by $\mO_d$ the Cuntz algebra generated by a $d$-dimensional Hilbert space $H$, and by $(\hrs)$ the vector space of linear operators from the tensor powers $H^r$ into $H^s$, $r,s \in \bN$ (for $r = 0$, we define $\iota :=$ $H^0 := \bC$). If $K \subseteq \ud$ is any closed subgroup, we define
\[
(\hrs)_K :=
\left\{
t \in (\hrs) : 
g_s t = t g_r
\ , \
g \in K
\right\}
\]
where $g_r := \otimes^r g \in (H^r , H^r)$, $g \in K$. We denote by $\wa K$ the category with objects $H^r$, $r \in \bN$, and arrows $(\hrs)_K$. It is well-known that $\wa K$ is a symmetric tensor \sC category, with symmetry induced by the flip operator $\theta \in (H^2 , H^2)_K$.

We also consider the well-known action
\begin{equation}
\label{def_ac}
\ud \ra {\bf aut} \mO_d
\ \ , \ \
u \mapsto \wa u : \wa u (t) := u_s t u_r^* 
\ , \ t \in (\hrs) \ .
\end{equation}

Let us denote by $\mO_K \subseteq \mO_d$ the \sC algebra generated by the set $\left\{ (\hrs)_K , \right.$ $\left. r,s \in \bN \right\}$. Then, it is easily verified that $\mO_K$ is the fixed-point algebra of $\mO_d$ w.r.t. the action (\ref{def_ac}) restricted to elements of $K$.

Let $\left\{ \psi_i \right\}_{i=1}^d$ be an orthonormal basis of $H$. If we consider the canonical endomorphism
\[
\sigma_d (t) := \sum_{i=1}^d \psi_i t \psi_i^* \ \ ,  \ \ t \in \mO_d \ ,
\]
then we find that $\sigma_d$ restricts to an endomorphism $\sigma_K \in {\bf end} \mO_K$.

For future reference, we introduce the notation
\[
{\bf aut}_{\sigma_K , \theta} \mO_K
:=
\left\{
\alpha \in {\bf aut}\mO_K
\ : \
\alpha \circ \sigma_K = \sigma_K \circ \alpha
\ , \
\alpha (\theta) = \theta
\right\}
\ \ .
\]

Let us now consider a normalized vector generating the totally 
antisymmetric tensor power $\bigwedge^d H$, say $R$. 
Then $R$ appears in the Cuntz algebra $\mO_d$ as an isometry 
$R \in$ $( \iota , H^d ) =$ $( \iota , \sigma_d^d )$, 
with support the totally antisymmetric projection $P_{\theta,d} :=$ 
$\sum_p sign (p) \theta (p)$, and satisfying the special conjugate equation
\begin{equation}
\label{eq_sce}
R^* \sigma_d (R) = (-1)^{d-1} d^{-1} 1
\end{equation}
Since $g_d R =$ $\det g R$, $g \in \ud$, we find that if $K \subseteq \sud$ then $R \in \mO_K$.

\

With $K \subseteq \ud$ as above, we define $NK$ as the normalizer of $K$ in $\ud$, and $QK$ as the quotient $NK / K$; we also denote by $p : NK \ra QK$ the natural projection. It is clear that $NK$ and $QK$ are compact Lie groups.

We denote by 
\[
{\bf aut} ( \mO_d , \mO_K )
\]
the group of automorphisms of $\mO_d$ leaving $\mO_K$ globally stable, and coinciding with the identity on $\mO_\ud \subseteq \mO_K$.

\begin{lem}
\label{lem_cd1}
Let $K \subseteq \sud$ be a closed group. Then, the following properties are satisfied:
\begin{enumerate}
\item $K$ is isomorphic to the stabilizer of $\mO_K$ in $\mO_d$, via the map (\ref{def_ac});
\item there is a commutative diagram
      \begin{equation}
      \xymatrix{
           {\bf aut} ( \mO_d , \mO_K )
      	    \ar[r]^-{\simeq}
      	    \ar[d]^-{\pi}
      	 &  NK
      	    \ar[d]^-{p}
      	 \\ {\bf aut}_{\sigma_K ,\theta} \mO_K
      	    \ar[r]^-{\simeq}
      	 &  QK
      }
      \end{equation}
      where the horizontal arrows are group isomorphisms, and the vertical arrows are group 
      epimorhisms.
\end{enumerate}
\end{lem}

\begin{proof}
See \cite{DR87}, \cite[Thm.5.14]{Vas3}.
\end{proof}
The element of ${\bf aut}_{\sigma_K ,\theta} \mO_K$ associated with $y \in QK$ will be denoted by $\wa y$.

\begin{rem}
Here we somehow provide a link to the situation discussed in \cite{CoDA}.
An element $\wa y \in {\bf aut}_{\sigma_K ,\theta} \mO_K$ extends to a $K$-{\em equivariant} automorphisms in ${\bf aut} ( \mO_d , \mO_K )$ if for some (and thus, for all) $u \in p^{-1}(y) \in NK$ the corresponding automorphism of $K$,
$\alpha_u := u \cdot u^*$, is inner.
In fact, in that case one readily finds that $u = g_0 u_0$ where $g_0 \in K$ and $u_0 \in CK$, the centralizer of $K$ in $\ud$.
In particular, if the only automorphisms of $K$ that are unitarily implemented in the defining representation $U$ of $K$ on $H$ are the inner ones, every element in ${\bf aut}_{\sigma_K ,\theta} \mO_K$ lifts to a $K$-equivariant automorphisms in ${\bf aut} ( \mO_d , \mO_K )$. (To see this, just notice that for $u \in NK$ one has $ug = (ugu^*)u$ for all $g \in K$, that is $u \in (U , U \circ \alpha_u)$.)

Perhaps the simplest example is provided by $K = \sud$ with $d=2$; for it is easy to see that the only equivariant extensions are given by the gauge automorphisms of $\mO_2$.
%
%
%
%
%
In turn,
the same conclusion holds true for 
the fundamental
representation of $\sud$, with $d \geq 2$
and, more generally,
for any closed subgroup $K \subseteq \sud$ that is quasi-complete and acts irreducibly on $H$.
\end{rem}

Now, let $X$ be a compact space.
We consider the rank $d$, free Hilbert $C(X)$-bimodule $E :=$ $C(X) \otimes H$, endowed with the natural left and right actions; it is clear that we may regard $E$ as the space of continuous maps from $X$ into $H$. The group of unitary $C(X)$-module operators of $E$ is given by the set
\[
C( X , \ud )
\]
of continuous maps from $X$ into the unitary group $\ud$, endowed with the topology of uniform convergence. 

For every $r,s \in \bN$, we denote by $(\ers)$ the Banach $C(X)$-bimodule of operators from the internal tensor power $E^r$ into $E^s$ (we also set $E^0 := C(X)$). It is clear that 
\[
( \ers ) \simeq C(X) \otimes (\hrs) \ \ .
\]
If $G$ is a closed subgroup of the unitary group $UE \equiv C( X , \ud )$, then in the same way as above we can consider the Banach $C(X)$-bimodules
\[
(\ers)_G := 
\left\{
t \in (\ers) : g_s t = t g_r
\ , \
g \in G
\right\}
\]
and denote by $\wa G_E$ the tensor \sC category with objects $E^r$, $r \in \bN$, and arrows $( \ers )_G$. If we identify $\theta \in (H^2 , H^2)$ with the corresponding constant map in $C(X) \otimes (H^2 , H^2)$, then we find that $\wa G_E$ is symmetric.

We denote by $\cog$ the \sC subalgebra of the Cuntz-Pimsner algebra $\coe \simeq C(X) \otimes \mO_d$ generated by the spaces $\left\{ (\ers)_G , r,s \in \bN \right\}$.
The canonical endomorphism $\sigma_d \in {\bf end} \mO_d$ naturally extends to give an endomorphism
\[
\sigma_E := id_X \otimes \sigma_d \in {\bf end} \coe
\]
(where $id_X \in$ ${\bf aut} C(X)$ denotes the identity automorphism); moreover, $\sigma_E$ restricts in a natural way to an endomorphism
\[
\sigma_G \in {\bf end} \cog \ \ .
\]
The \sC dynamical systems $( \coe , \sigma_E )$, $( \cog , \sigma_G )$ have been studied in a more general setting in \cite{Vas}. There is an action
\begin{equation}
\label{def_ace}
UE \ra {\bf aut} \coe 
\ \ , \ \
u \mapsto \wa u 
\ : \ 
\wa u (t) := u_s t u_r^* \ , \ t \in (\ers) 
\ \ .
\end{equation}
If we restricts (\ref{def_ace}) to elements of $G$, then by construction we find that $\cog$ is contained in the fixed-point algebra of $\coe$ w.r.t. this restricted action. 


\begin{rem}
In general, the space $(\ers)_G \subset$ $C(X) \otimes (\hrs)$ {\em is not} a free $C(X)$-bimodule. For example, in the case in which $G$ is a locally compact group, then the Haar measure induces a projection
\[
P_G : C(X) \otimes (\hrs) \to (\ers)_G
\ ,
\]
and this implies that in general $(\ers)_G$ is projective (and not free). 
We give an explicit example. Let us consider the $2$-sphere $S^2$, and let $\mathcal E \to S^2$ denote the (complexified) tangent bundle. It is well-known that $\mathcal E$ is nontrivial (i.e., the module of sections of $\mathcal E$ is not a free $C(S^2)$-module); on the other side, if $V^n := S^2 \times \bC^n$, $n \in \bN$, denotes the trivial rank $n$ vector bundle, then we find $\mathcal E \oplus V^1 \simeq V^3$ (see \cite[II.1.19]{Kar}). 
Let us now denote by $\bT$ the torus, acting in the natural way on $V^1$. If we consider the group $G :=$ $1 \oplus \bT$, where $1$ is the identity on $\mathcal \mE$, then there is a natural action
\[
G \times (\mathcal E \oplus V^1) 
\ \to \ 
\mathcal E \oplus V^1
\ \ ,
\]
such that the fixed-point vector bundle is $\mathcal E$. By passing to the modules of sections, this means that if we define $M :=$ $C(S^2) \otimes \bC^3$, then we obtain that there is a $G$-action
\[
G \times M \to M
\]
such that $( \iota , M )_G$ is not free as a $C(S^2)$-module (in fact, $( \iota , M )_G$ is isomorphic to the module of sections of $\mE$).
\end{rem}

It is clear that $\cog$ is a continuous field of \sC algebras with base space $X$. The above considerations imply that in general $\cog$ is not trivial as a continuous field of \sC algebras, anyway there is a $C(X)$-monomorphism $\cog \hra \coe =$ $C(X) \otimes \mO_d$.
We now analyze the structure of the fibres of $\cog$. At first, we note that for every $x \in X$ it turns out that the fibre $(\cog)_x :=$ $\cog / ( C_x(X) \cog )$ is naturally embedded in $\mO_d$. We define
\[
G^x :=
\left\{
u \in \ud : \wa u (t) = t \ \forall t \in (\cog)_x
\right\}
\ \ .
\]
Clearly, $G^x$ is a closed subgroup of $\ud$.
%
\begin{lem}
\label{lem_ogx}
For every $x \in X$, one has $(\cog)_x \simeq \mO_{G^x}$.
\end{lem}

\begin{proof}
See \cite[Cor.4.6]{Vas}.
\end{proof}

We now define the {\em spectral bundle} 
\[
\mG := \{(x,y) \ | \ x \in X, y \in G^x\} 
\subseteq X \times \ud
\ \ ,
\]
endowed with the obvious topology as a subspace of $X \times \ud$.
Moreover, we introduce the space of sections of $\mG$
\[
SG := 
\left\{
g \in C(X,\ud) : g(x) \in G^x, \ \forall x \in X \ 
\right\} .
\]
As in \cite[Lemma 4.10]{Vas}, we find
\begin{lem}
There is an inclusion $G \subseteq SG$; moreover, $(\ers)_G = (\ers)_{ SG}$, $r,s \in \bN$, 
and $\cog =$ $\mO_{SG}$.
\end{lem}

\begin{rem}
Since $SG$ may be recovered as the stabilizer in $\coe$ of $\cog = \mO_{SG}$ (see \cite[Prop.4.8]{Vas}), we find that $SG =$ $S (SG)$ (in fact, $\mO_{SG} = \mO_{S(SG)}$).
\end{rem}

Let $NG^x$ denote the normalizer of $G^x$ in $\ud$. We define
\[
\mNG
:= 
\{(x,n) \ | \ x \in X, n \in NG^x \} \subseteq X \times \ud
\]
endowed with the natural topology as a subspace of $X \times \ud$, and the quotient space (with the associated natural projection)
\[
\mQG := \mNG / \mG 
\ \ , \ \ 
p : \mNG \ra \mQG
\ \ .
\]
It is clear that, at the level of sets, there is an identification $\mQG =$ 
$\{(x,z) \ | \ x \in X, z \in NG^x / G^x \}$.
We define $q : \mQG \ra X$ as the natural projection. For every closed group $G \subseteq C(X,\ud)$, we define
\[
NG :=
\left\{
u \in C(X,\ud) :
u(x) \in NG^x, \ \forall x \in X
\right\} \ . 
\] 
It is clear that $G$ is a normal subgroup of $NG$. We also define the set of continuous sections of $\mQG$
\[
QG := 
\left\{
y : X \hra \mQG : q \circ y = id_X
\right\}
\ \ .
\]
Note that there is a morphism
\begin{equation}
\label{def_pnq}
p_* : NG \ra QG \ \ , \ \ p_*(u) := p \circ u \ \ .
\end{equation}
In general, $p_*$ is not surjective. 
For a short discussion of this point
in the simple
but instructive
case where $G = C(X,K)$ for some closed subgroup $K \subset \sud$, 
see the Appendix.

\

We now give a version of Lemma \ref{lem_cd1} for the \sC algebra $\coe$. For this purpose, we introduce
\[
{\bf aut}_{\sigma_G , \theta} \cog
:=
\left\{
\alpha \in {\bf aut} \cog
\ : \
\alpha \circ \sigma_G = \sigma_G \circ \alpha
\ , \
\alpha (\theta) = \theta
\right\}
\ \ ,
\]
and
\[
{\bf aut} ( \coe , \cog )
:=
\left\{
\alpha \in {\bf aut} \coe :
\alpha (\cog) = \cog
, \ 
\alpha |_{\mO_{UE}} = id
\right\}
\ \ .
\]

\begin{lem}
\label{lem_cd2}
Let $G \subseteq C(X,\sud)$ be a closed group. Then, the following properties are satisfied:
\begin{enumerate}
\item $SG$ is isomorphic to the stabilizer of $\cog$ in $\coe$, 
via the map (\ref{def_ace});
\item there is a commutative diagram
      \begin{equation}
      \xymatrix{
           {\bf aut} ( \coe , \cog )
      	    \ar[r]^-{\simeq}
      	    \ar[d]^-{\pi}
      	 &  NG
      	    \ar[d]^-{p_*}
      	 \\ {\bf aut}_{\sigma_G ,\theta} \cog
      	    \ar[r]^-{\simeq}
      	 &  QG
      }
      \end{equation}
      where the horizontal arrows are group isomorphisms.
\end{enumerate}
\end{lem}

\begin{proof}
Point (1) is a consequence of \cite[Prop.4.8]{Vas}. Concerning Point (2), we note that every $\alpha \in$ ${\bf aut} ( \coe , \cog )$ is by definition a $C(X)$-automorphism of $\coe$, thus there exists a continuous family $\alpha_x \in {\bf aut} \mO_d$ such that $\pi_x \circ \alpha =$ $\alpha_x \circ \pi_x$, where $\pi_x : \coe \ra \mO_d$, $x \in X$, is the evaluation epimorphism (see \cite[\S 4]{Nil96}). Since $\pi_x \circ \sigma_E =$ $\sigma_d \circ \pi_x$, $\pi_x (\cog) =$ $\mO_{G^x} =$ $(\cog)_x$ (see Lemma \ref{lem_ogx}), we find that each $\alpha_x$ belongs to ${\bf aut} ( \mO_d , \mO_{G^x} )$. By Lemma \ref{lem_cd1}, we conclude that there is $u_x \in NG^x$ such that $\alpha_x = \wa u_x$. Since the family $(\alpha_x)_x$ is continuous, and from the fact that the correspondence (\ref{def_ac}) is one-to-one, we conclude that the family $(u_x)_x$ is unique and continuous. Thus, we have proved that 
\[
NG \ni u \mapsto \wa u \in {\bf aut} ( \coe , \cog )
\]
defines the desired automorphism. Now, if $v \in NG^x$, $y \in G^x$, we find $\wa{vy} \circ \pi_x (t) = \wa u \circ \pi_x (t)$ for every $t \in \cog$; in other terms, $NG^x$ acts on $\mO_{G^x}$. This implies that there is a well-defined map $QG \ra {\bf aut}_{\sigma_G ,\theta} \cog$, which assigns to the generic element $y \in QG$ the automorphism $\wa y$. If $\wa y = \wa {y'}$, then $\wa{ y^{-1}y' }$ is the identity on $\cog$, and this means that for every $x \in X$ there is $g_x \in G^x$ such that $y'(x) = y(x)g_x$; of course, by definition of $NG$ this implies $y = y'$. It remains to verify that $\left\{ y \mapsto \wa y \right\}$ is surjective. But this easily follows from the fact that every $\beta \in {\bf aut}_{\sigma_G ,\theta} \cog$ defines a continuous family $(\beta_x)_x$, $\beta_x \in {\bf aut}_{\sigma_K ,\theta} \mO_K$, and from Lemma \ref{lem_cd1}.
\end{proof}

\section{Crossed Products}

Let $\mA$ be a \sC algebra with centre $\mZ$, $\rho \in {\bf end} \mA$ a endomorphism. We consider the tensor \sC category $\wa \rho$ having as objects the powers $\rho^r$, $r \in \bN$ (for $r = 0$, we consider the identity automorphism $\iota := \rho^0$), and arrows the intertwiner spaces
\[
(\rhors) = \left\{  t \in \mA : \rho^s (a) t = t \rho^r (a) \ , \ a \in \mA  \right\} \ ;
\]
\noindent the tensor product is given by
\[
\left\{ 
\begin{array}{l}
\rho^r  ,  \rho^s    \mapsto     \rho^{r + s}  \\ 
t , s    \mapsto   t  \rho^r (t') = \rho^s (t') t \in (\rho^{r+r'} , \rho^{s+s'})
\end{array}
\right.
\]
\noindent $t \in (\rho^r,\rho^s) , t' \in (\rho^{r'},\rho^{s'})$. Moreover, we introduce the notation
\[
\mZ^\rho := \left\{ f \in \mZ : \rho (f) = f \right\}
\]
\noindent and denote by $X^\rho$ the spectrum of $\mZ^\rho$, so that we make the identification $\mZ^\rho \simeq \zro$. Note that $\mA$ has a natural structure of $\zro$-algebra in the sense of Kasparov (i.e., there is a nondegenerate morphism from $\zro$ into $\mZ$). In general, $\zro$ does not coincide with $\mZ$.

In order to give an intuitive idea of the role of $\zro$, we consider the case in which $\rho$ restricts to an automorphism of $\mZ$. In this case, by Gel'fand transform we obtain a homeomorphism $\alpha_\rho$ of the spectrum $X'$ of $\mZ$, and $X^\rho$ coincides with the orbit space $X' / \bZ$ w.r.t. the $\bZ$-action induced by $\alpha_\rho$.

\begin{defn}[Permutation symmetry]
\label{def11}
An endomorphism $\rho$ of $\mA$ has permutation symmetry (PS) if there is a unitary representation $p \mapsto \eps (p)$ of the group $\bP_\infty$ of finite permutations of $\bN$ in $\mA$, such that:
\begin{equation}\label{ps1} \eps (\bS p) = \rho \circ \eps (p) \end{equation}
\begin{equation}\label{ps2} \eps := \eps (1,1) \in (\rho^2 , \rho^2) \end{equation}
\begin{equation}\label{ps3} \eps (s,1) \ t = \rho (t) \ \eps (r,1) \ \ , \ \ t \in (\rhors) \ , \end{equation}
\noindent where $(r,s) \in \bP_{r+s}$ permutes the first $r$ terms with the remaining $s$, and $\bS$ is the shift $(\bS p)(1) := 1$, $(\bS p)(n) := 1 + p(n-1)$, $p \in \bP_\infty$. The above properties imply that
\begin{equation}
\label{weak_perm}
\eps(p) \in (\rho^n , \rho^n) \ \ , \ \ n \in \bN \ , \ p \in \bP_n \ ;
\end{equation}
\noindent we say that $\rho$ has {\bf weak permutation symmetry} (WPS) if just (\ref{ps1}),(\ref{weak_perm}) hold. 
\end{defn}


\begin{defn}
\label{def_strict_int}
Let $\rho$ be an endomorphism on a \sC algebra $\mA$ carrying a weak permutation symmetry $\eps$. The elements of $(\rhors)$ for which (\ref{ps3}) holds are said {\bf symmetry intertwiners}:
\[
(\rhors)_\eps := 
\left\{ t \in (\rhors) : \eps(s,1) t = \rho(t) \eps (r,1) \right\} \ .
\]
\end{defn}
It is a trivial consequence of the definition that
$\eps(p) \in (\rho^n,\rho^n)_\eps$ for all $p \in \bP_n$, if $\rho$ has a WPS. 
Also, $(\rhors)_\eps = (\rhors)$, $r,s \in \bN$ if the WPS is a PS. 

\begin{rem}
We denote by $\wa \rho_\eps$ the \sC subcategory of $\wa \rho$ having arrows $(\rhors)_\eps$, $r,s \in \bN$. It is easily verified that $\wa \rho_\eps$ is a tensor subcategory of $\wa \rho$ (see \cite[Lemma 4.2]{Vas2}). In particular, by definition $(\ii)_\eps = \zro$, so that every $(\rhors)_\eps$ has a natural structure of Banach $\zro$-bimodule w.r.t. the multiplication by elements of $\zro$.
\end{rem}

\begin{defn}[Permutation quasi-symmetry]
\label{def_gps}
An endomorphism $\rho$ of a \sC algebra $\mA$ with centre $\mZ$ has permutation quasi-symmetry (qPS) if it has WPS and 
\begin{equation}
\label{gps3} 
(\rhors) = \rho^s (\mZ) \cdot (\rhors)_\eps = (\rhors)_\eps \cdot \rho^r (\mZ) 
\ . 
\end{equation}
\end{defn}

The previous equation has to be intended in the sense that the set $\left\{ zt , z \in \mZ , \right.$ $\left. t \in (\rhors)_\eps \right\}$ is dense in $(\rhors)$, $r,s \in \bN$.

\begin{ex}
The endomorphisms $\sigma_d \in {\bf end} \mO_d$, $\sigma_E \in {\bf end} \coe$, $\sigma_G \in {\bf end} \cog$ considered in Sec.\ref{intro} have PS.
\end{ex}

\begin{ex}
The so-called {\em canonical endomorphisms} considered in \cite{BL} have qPS. In that particular case, every $(\rhors)_\eps$ is a free $\zro$-bimodule.
\end{ex}

Let $\rho \in {\bf end} \mA$ be an endomorphism with WPS. Then $\rho$ has a well-defined {\em dimension} $d(\rho) \in \bN$ (see \cite{Vas,DR87}). The data of an endomorphism with WPS $\eps$ and dimension $d$ will be denoted by $( \rho , \eps , d )$.

We now introduce the following \sC subalgebras of $\mA$:
\begin{equation}
\label{def_oro}
\oro := C^* \left\{ t , t \in (\rhors) , r,s \in \bN \right\}
\ ;
\end{equation}
\begin{equation}
\label{def_oroeps}
\oroeps := C^* \left\{ t , t \in (\rhors)_\eps , r,s \in \bN \right\} 
\ ;
\end{equation}
\begin{equation}
\label{def_aeps}
\mP_{\rho,\eps} := C^* \left\{  \eps (p) , p \in \bP_\infty \right\} 
\ .
\end{equation}
Then clearly $\mP_{\rho,\eps} \subseteq \oroeps \subseteq \oro \subseteq \mA$.
It is proved in \cite{DR87} that there is an isomorphism
\begin{equation}
\label{def_i}
i : \mO_\ud \ra \mP_{\rho,\eps} \ \ .
\end{equation}
that is, in essence, just Weyl reciprocity.

\begin{rem}
\label{rem_WQPS}
Of course, the implications $PS$ $\Rightarrow$ $qPS$ $\Rightarrow$ $WPS$ hold. The above-defined 'soft' notions of permutation symmetry have been considered in \cite{Vas2}, and the interest in them arises from the fact that several properties can be proved by assuming WPS or qPS. For example, it suffices to assume WPS to prove that $\oroeps$ is a continuous bundle with fibres \sC algebras of the type $\mO_K$, $K \subseteq \sud$ (\cite[Thm.5.1]{Vas2}).
Moreover, it suffices to assume qPS to prove that the intertwiner spaces $(\rhors)$, $r,s \in \bN$, can be interpreted in terms of equivariant operators between tensor powers of a suitable Hilbert $\mZ$-bimodule (\cite[Thm.7.4]{Vas0}). 
\end{rem}

We introduce the following notation. If $\rho' \in {\bf end} \mA'$ has WPS, say $\eps'$, then a \sC algebra morphism $\eta : \mA \ra \mA'$ is said {\em equivariant} if
\[
\rho' \circ \eta = \eta \circ \rho \ \ , \ \ \eta (\eps) = \eps' \ \ .
\]
In such a case, we use the notation $\eta : ( \mA , \rho , \eps ) \ra ( \mA' , \rho' , \eps' )$. The group of equivariant automorphisms $\alpha : ( \mA , \rho , \eps ) \ra ( \mA , \rho , \eps )$ will be denoted by
\begin{equation}
\label{def_autro}
{\bf aut}_{\rho,\eps} \mA \ \ .
\end{equation}
Then ${\bf aut}_{\rho,\eps} \mA$ is a (not normal, in general) subgroup of ${\bf aut} \mA$.

\begin{rem}
\label{rem_autro}
Let $\alpha \in {\bf aut}_{\rho, \eps} \mA$. Then,
\begin{enumerate}
\item $\alpha (\rhors) = (\rhors)$, $r,s \in \bN$;
\item $\alpha (\rhors)_\eps = (\rhors)_\eps$, $r,s \in \bN$;
\item $\alpha (\eps(p)) = \eps (p)$, $p \in \bP_\infty$.
\end{enumerate}
\end{rem}

\

Let $\mA$ be a \sC algebra with centre $\mZ$, and $\rho$ an endomorphism with WPS, say $\eps$, with dimension $d := d(\rho)$. A {\em minimal crossed product} of $\mA$ by $\rho$ (mCP) is given by a \sC algebra $\mB$ with identity $1 \in \mB$, satisfying the following properties:
\begin{enumerate}
\item  $\mB$ is generated by $\mA$ and a set $\left\{ \psi_i \right\}_{i=1}^d$ of $d$ isometries, satisfying the Cuntz relations
\begin{equation}
\label{eq_cuntz}
\psi_i^* \psi_j = \delta_{ij} 1 \ \ , \ \ \sum_i \psi_i \psi_i^* = 1 \ \ .
\end{equation}
By universality of the Cuntz algebra, this implies that there is a monomorphism $j : \mO_d \hra \mB$. Moreover, the following endomorphism of $\mB$ is defined:
\begin{equation}
\label{def_sb}
\sigma (b) := \sum_i \psi_i b \psi_i^* \ \ ;
\end{equation}
\item $\sigma (a) = \rho (a)$, $a \in \mA$;
\item $\mA' \cap \mB = \mZ$ (minimality);
\item  $i (t) = j(t)$, $t \in \mO_\ud \subset \mO_d$, where $i$ is defined in (\ref{def_i}) (symmetry).
\end{enumerate}
Note that there is a naturally defined WPS $( \sigma , \eps , d )$, thus it makes sense to consider the group 
\[
{\bf aut}_{\sigma,\eps} \mB \ \ .
\]
Let us consider the group $C( X^\rho , \ud )$ (see Sec.\ref{intro}). If $E = C(X^\rho,H)$ is the (free) Hilbert $\zro$-bimodule introduced in Sec.\ref{intro}, then it is clear that $C( X^\rho , \ud )$ coincides with the unitary group of $E$. In the sequel, we will identify $E$ with
\[
{\mathrm{span}} \ \left\{ \psi_i f , \ i = 1 , \ldots , d , \ f \in \zro \right\} 
\]
so that there is an inclusion $E \subset \mB$. The \sC subalgebra of $\mB$ generated by $E$ is clearly isomorphic to $\coe \simeq \zro \otimes \mO_d$; in this way, $j$ extends to a $\zro$-monomorphism
\begin{equation}
\label{eq_j}
j : ( \coe , \sigma_E , \theta ) \hra ( \mB , \sigma , \eps )
\ .
\end{equation}
Further, one readily checks that $j \circ \sigma_E = \sigma \circ j$ and $j(\theta) = \epsilon$.
\begin{lem}
\label{lem_eib}
With the above notation, it turns out $E =$ $( \iota_\mB , \sigma )$, where $\iota_\mB \in$ ${\bf aut} \mB$ denotes the identity automorphism.
\end{lem}
\begin{proof}
By (\ref{def_sb}) it follows that $E \subseteq$ $( \iota_\mB , \sigma )$. On the converse, let $b \in$ $( \iota_\mB , \sigma )$; then, for every $i = 1 , \ldots , d$ we define $b_i := \psi_i^* b$. In order to prove that $b \in E$, it suffices to verify that $b_i \in \zro$. For this purpose, we note that 
\[
b' b_i = b' \psi_i^* b = \psi_i^* \sigma (b') b = b_i b'
\ ,
\]
so that $b_i \in \mB' \cap \mB$. Since $\mB' \cap \mB \subseteq$ $\mA' \cap \mB$, the minimality condition implies that $b_i \in \mZ$. Moreover, we have
\[
\rho (b_i) = 
\sigma (b_i) = 
\sum_k \psi_k b_i \psi_k^* = 
b_i
\ ,
\]
and this implies $b_i \in \zro$.
\end{proof}

\begin{thm}
\label{thm_dual}
Let $\mA$ be a \sC algebra with centre $\mZ$, and $( \rho , \eps , d )$ a quasi-symmetric endomorphism. For every mCP $\mB$, there exist:
\begin{enumerate}
\item A closed subgrop $G$ of $C( X^\rho , \ud )$, naturally acting on the $\zro$-bimodule $E$.
\item An equivariant monomorphism 
\begin{equation}
\label{eq_mu}
\mu : ( \cog , \sigma_G , \theta ) \to ( \mA , \rho , \eps  ) \ \ ;
\end{equation}
\item An isomorphism $( \wa G_E , \theta ) \to ( \wa \rho_\eps , \eps )$ of symmetric tensor \sC categories.
\end{enumerate}
\end{thm}
\begin{proof}
{\bf Point 1}: We define 
\[
G :=
{\bf aut}_{\mA,\rho} \mB
:=
\left\{
\beta \in {\bf aut} \mB : \beta |_\mA = id , \ \beta \circ \sigma = \sigma \circ \beta
\right\} \ ,
\]
and show that $G$ acts in a natural way on $E$. This will suffice to conclude that there is an inclusion $G \hra C( X^\rho , \ud )$. Let $\psi \in E =$ $( \iota_\mB , \sigma )$ (recall Lemma \ref{lem_eib}), and $\beta \in {\bf aut}_{\mA,\rho} \mB$. For every $i = 1 , \ldots , d$, we define $b_i := \psi_i^* \beta (\psi)$, and note that for every $a \in \mA$ it turns out
\[
a b_i = 
\psi_i^* \rho (a) \beta (\psi) =
\psi_i^* \beta ( \rho (a) \psi ) =
\psi_i^* \beta (\psi) a =
b_i a
\ .
\]
This implies that $b_i \in $ $\mA' \cap \mB = \mZ$. Moreover, by using $\sigma (\psi) = \eps \psi$ and $\beta (\eps) = \eps$ we find
\[
\rho (b_i) =
\sigma(\psi_i^*) \beta ( \sigma(\psi) ) =
\psi_i^* \eps \beta ( \eps \psi ) =
\psi_i^* \eps^2 \beta (\psi) =
b_i
\ ,
\]
and this implies $b_i \in \zro$. Thus, $\beta (\psi) = \sum_i \psi_i b_i$ belongs to $E$, and this implies that $\beta |_E$ restricts to a unitary operator on $E$. By applying (\ref{def_ace}), we find that there is $u \in C(X^\rho , \ud )$ such that 
\begin{equation}
\label{eq_ub}
\beta \circ j = j \circ \wa u
\ \ , \ \ 
\wa u \in {\bf aut} \coe
\ .
\end{equation}
This proves Point 1.\\
{\bf Point 2}: Let $\cog \subseteq \coe$ be the \sC algebra generated by the spaces $(\ers)_G$. We define $\mu := j |_{\cog}$, $\mu : \cog \hra \mB$. By (\ref{eq_j}), $\mu$ is equivariant:
\[
\mu : ( \cog , \sigma_G , \theta ) \hra ( \mB , \sigma , \eps )
\ .
\]
Moreover, by applying (\ref{eq_ub}) we find $\mu (\ers)_G =$ $\mu ( \sigma_G^r , \sigma_G^s ) \subseteq$ $(\rhors)_\eps$, $r,s \in \bN$. 
In order to prove the opposite inclusion, let us define
\begin{equation}
\label{def_psiI}
\psi_I := \psi_{i_1} \cdots \psi_{i_r} \in 
E^r = j( \iota , \sigma^r ) \ .
\end{equation}
It is clear that
\[
\psi_I^* \psi_J = \delta_{IJ} 1
\ \ , \ \ 
\sum_I \psi_I \psi_I^* = 1
\ .
\]
If $t \in (\rhors)_\eps$, then we find $t = \sum_{IJ} \psi_I t_{IJ} \psi_J^*$, where $t_{IJ} := \psi_I^* t \psi_J$. Since $\sigma |_\mA = \rho$, we find 
\[
t_{IJ} a = 
\psi_I^* t \rho^r(a) \psi_J = 
\psi_I^* \rho^s(a) t \psi_J =
a t_{IJ}
\ \ ,
\]
thus $t_{IJ} \in \mA' \cap \mB = \mZ$. Moreover, since $E^r \subseteq ( \iota_\mB , \sigma^r )_\eps$, by applying Def.\ref{def_strict_int} we find
\[
\rho (t_{IJ}) =
\psi_I^* \eps (1,s) \eps (s,1) t \eps (1,r) \eps (r,1) \psi_J =
t_{IJ}
\ \ .
\]
This implies $t_{IJ} \in \zro$, thus $t \in j (\ers)$. Moreover, (\ref{eq_ub}) implies that $\wa u \circ j^{-1} (t) = j^{-1} (t)$ for every $u \in G$. Thus, we conclude that $(\rhors)_\eps =$ $j( \ers )_G$, $r,s \in \bN$.\\
{\bf Point 3} follows trivially from Point 2.
\end{proof}

\begin{rem}\label{cogoroeps}
In particular, it follows from the last proof that $\mu(\cog) = \oroeps$.
\end{rem}

\begin{rem}
If $\mZ \neq \bC 1$, then existence and unicity of $( \mB , G )$ are not ensured. 
For a class of examples of this phenomenon, see \cite[\S 6.2]{Vas0}.
In the particular case studied in \cite{BL}, then $( \mB , G )$ exists and is unique, with the additional property that $G$ is of the type $G =$ $C(X^\rho ,K)$ for some $K \subseteq \ud$. For a proof of this result, see \cite[Thm.7.4]{Vas0}.
\end{rem}

Let $( \rho , \eps , d )$ be an endomorphism with WPS (resp. qPS, PS). We say that $\rho$ is {\em weakly special} (resp. {\em quasi special}, {\em special}) if there exists an isometry $S \in$ $( \iota , \rho^d)$ such that
\begin{equation}
\label{eq_scp}
RR^* = P_{\eps,d} := \sum_{p \in \bP_d} sign (p) \ \eps (p)
\ \ , \ \
R^* \rho (R) = (-1)^{d-1} d^{-1} 1 \ \ .
\end{equation}
The data of a weakly special endomorphism will be denoted by $( \rho , \eps , d , R )$.

\begin{ex}
If $K \subseteq \sud$, $G \subseteq C (X,\sud)$, then the endomorphisms considered in Sec.\ref{intro} are special (see (\ref{eq_sce})). Special endomorphisms in \sC algebras with trivial centre have been studied in \cite{DR89A}. Moreover, if a canonical endomorphism in the sense of  \cite{BL} satifies (\ref{eq_scp}), then it is quasi-special.
\end{ex}

\begin{prop}
Let $( \rho , \eps , d , R )$ be a canonical endomorphism of a \sC algebra $\mA$ (in the sense of
\cite{BL}) satisfying (\ref{eq_scp}). Then, there exists a mCP $\mB$ of $\mA$ by $\rho$, with $G = C(X,K)$ for some closed subgroup $K \subseteq \sud$.
\end{prop}

\begin{proof}
See \cite[Thm.7.4]{Vas0}.
\end{proof}

\section{Extensions}

Let $( \rho , \eps , d , S )$ be a quasi-special endomorphism of a \sC algebra $\mA$. If $E$ is the free, rank $d$ Hilbert $\zro$-bimodule, then we refer the reader to the notations $NG$, $QG$ introduced in Sec.\ref{intro} for every closed group $G \subseteq C(X^\rho,\ud)$.

Let $\mB$ be a mCP of $\mA$ by $\rho$. We recall that $\mB$ is endowed with a weakly special endomorphism $( \sigma , \eps , d , S )$ extending $\rho$, and denote by
\[
{\bf aut}_{\sigma,\eps} ( \mB ; \mA )
\]
the group of automorphisms of $\mB$ commuting with $\sigma$, leaving $\eps$ fixed, and $\mA$ globally stable.

Now, by (\ref{eq_j}) and Thm.\ref{thm_dual}, there is a commutative diagram
\begin{equation}
\label{30}
\xymatrix{
            \mA
	    \ar[r]^{\subset}
	 &  \mB
	 \\ \cog
	    \ar[u]^{\mu}
	    \ar[r]^{\subset}
	 &  \coe
	    \ar[u]_{j}
}
\end{equation}
The above diagram, together with the condition $\rho (a) = \sigma (a)$, $a \in \mA$ (and (\ref{def_sb})), implies that $\mB$ may be regarded as the crossed product $\mA \rtimes_\mu \wa G$ in the sense of \cite[\S 3]{Vas2}. 
Universality property of the crossed product implies that if $\mB'$ is a \sC algebra satisfying (\ref{30}) and (\ref{def_sb}) with given \sC monomorphisms $\mu' : \cog \to \mA$, $j' : \coe \to \mB'$, then there is an isomorphism 
\begin{equation}
\label{eq_beta}
\beta : \mB \to \mB' \ ,
\end{equation}
such that $\beta \circ j =$ $j' \circ \beta$ (this also implies $\beta \circ \mu =$ $\mu' \circ \beta$); moreover, $\beta$ turns out to be also an $\mA$-module map, i.e. $\beta (\mA) = \mA$.
The previous universal property is the key ingredient of the next result.

\begin{thm}
Let $\mA$ be a \sC algebra, $( \rho , \eps , d , S)$ a quasi-special endomorphism, and $\mB$ a mCP of $\mA$ by $\rho$ with a group $G \subseteq$ $C(X^\rho,\sud)$, $G \simeq$ ${\bf aut}_{\mA,\rho}\mB$. Then, there is a commutative diagram
\begin{equation}
\label{diag_aut}
\xymatrix{
            {\bf aut}_{\sigma,\eps} ( \mB ; \mA )
	    \ar[r]^-{r_\mB}
	    \ar[d]^-{\pi}
	 &  NG
	    \ar[d]^-{p_*}
	 \\ {\bf aut}_{\rho,\eps} \mA
	    \ar[r]^-{r_\mA}
	 &  QG
}
\end{equation}
In addition, $r_\mB$ is a group epimorphism. 
\end{thm}

\begin{proof}
The proof of the theorem will be divided in several steps.

{\bf Step 1.} Let $\beta \in {\bf aut}_{\sigma,\eps} ( \mB ; \mA )$. For every $t \in ( \sigma^r , \sigma^s )$, $b \in \mB$, we find
\[
\beta (t) \sigma^r (b) = 
\beta ( t \sigma^r \circ \beta^{-1} (b) ) =
\beta ( \sigma^s \circ \beta^{-1} (b) t ) =
\sigma^s (b) \beta (t) \ ,
\]
thus $\beta ( \sigma^r , \sigma^s ) =$ $( \sigma^r , \sigma^s )$.
In particular, it turns out $\beta ( \iota_\mB , \sigma ) = $ $( \iota_\mB , \sigma )$, thus by Lemma \ref{lem_eib} we conclude that $\beta$ restricts to a unitary map $u : E \to E$. Let us denote by $\wa u \in$ ${\bf aut} \coe$ the automorphism obtained by extending $u$ in the usual way. By (\ref{eq_j}), it turns out that 
\begin{equation}
\label{eq_bj}
\beta \circ j = j \circ \wa u \ .
\end{equation}

{\bf Step 2.} Our task is now to give an explicit representation of $\beta$ in terms of its restriction to $\mA$ and $\wa u$. As a first step, we verify that $\beta ( \rhors )_\eps =$ $( \rhors )_\eps$.
For this purpose, note that the condition $\beta ( \sigma^r , \sigma^s ) =$ $( \sigma^r , \sigma^s )$ is equivalent to $\beta \circ j (\ers) =$ $j ( \ers )$, $r,s \in \bN$ (in fact, Lemma \ref{lem_eib} implies $( \sigma^r , \sigma^s ) =$ $j (\ers)$). In particular, if $t \in (\rhors)_\eps$, then by Thm.\ref{thm_dual} we find $t \in \mu (\ers)_G \subseteq$ $j (\ers)$. Thus, $\beta(t) \in$ $j(\ers) =$ $(\sigma^r,\sigma^s)$, and we obtain the inclusion
\[
\beta (\rhors)_\eps \subseteq ( \rhors ) \ .
\]
Moreover, by (\ref{eq_j}), we find
\[
\rho (t) = 
\rho \circ j (t') =
\mu \circ \sigma_E (t') =
j ( \theta (s,1) t' \theta (1,r) )  =
\eps (s,1) t \eps (r,1)
\ .
\]
Thus $\beta (\rhors)_\eps \subseteq (\rhors)_\eps$, and by applying the above argument to $\beta^{-1}$ we conclude 
\begin{equation}
\label{eq_brse}
\beta (\rhors)_\eps = (\rhors)_\eps \ .
\end{equation}
Let us denote by $G \subseteq UE$ the closed group constructed in Thm.\ref{thm_dual}; by construction, $G$ satisfies the property that $j$ restricts to the isomorphism (\ref{eq_mu}), in such a way that $\mu (\ers)_G =$ $( \rhors )_\eps$, $r,s \in \bN$. Since $\beta ( \rhors )_\eps =$ $( \rhors )_\eps$ and $\beta \circ \mu =$ $\beta \circ j |_{\cog} =$ $j \circ \wa u$, we obtain
\begin{equation}
\label{eq_ueg}
\wa u ( \ers )_G = ( \ers )_G 
\ \ , \ \
r,s \in \bN 
\ .
\end{equation}
Thus, from Lemma \ref{lem_cd2}, we conclude that in particular $u \in NG$, and $\wa u$ restricts to an automorphism of $\cog$. 

{\bf Step 3.} The relations (\ref{eq_brse}) imply that by restricting to $\mA$ an element of ${\bf aut}_{\sigma,\eps} ( \mB ; \mA )$, one obtains an automorphism of $\mA$ preserving the spaces of symmetry intertwiners. Thus, the following map is well-defined:
\[
\pi : {\bf aut}_{\sigma,\eps} ( \mB ; \mA ) \to {\bf aut}_{\rho,\eps} \mA
\ \ , \ \
\pi (\beta) := \beta |_\mA
\ .
\]
Moreover, if $u \in NG$ is the unitary constructed as in the previous step, we define
\begin{equation}
\label{eq_ho1}
r_\mB : {\bf aut}_{\sigma,\eps} ( \mB ; \mA ) \to NG
\ \ , \ \ 
r_\mB (\beta) := u
\ .
\end{equation}
By (\ref{eq_bj}) and Lemma \ref{lem_cd2}, we find that $r_\mB$ preserves the product, thus it is a group morphism. Moreover, for every $u \in NG$ we have a commutative diagram
\begin{equation}
\label{31}
\xymatrix{
            \mA
	    \ar[r]^{\subset}
	 &  \mB
	 \\ \cog
	    \ar[u]^{\mu \circ \wa u }
	    \ar[r]^{\subset}
	 &  \coe
	    \ar[u]_{j \circ \wa u}
}
\end{equation}
By (\ref{eq_beta}), we conclude that there is an automorphism $\beta \in {\bf aut} \mB$ satisfying (\ref{eq_bj}); this means that $r_\mB$ is surjective.

{\bf Step 4.} Let $\alpha \in {\bf aut}_{\rho,\eps} \mA$. From the isomorphism (\ref{eq_mu}) 
and Remark \ref{cogoroeps},
we obtain that there is $\alpha' \in {\bf aut}_{\sigma_G , \theta} \cog$ such that $\alpha \circ \mu =$ $\mu \circ \alpha'$. By Lemma \ref{lem_cd2}, we find that there is $y \in QG$ such that $\alpha' = \wa y$. Thus, we can define define the morphism
\[
r_\mA : {\bf aut}_{\rho,\eps} \mA \to QG
\ \ , \ \
r_\mA (\alpha) := y
\ .
\]

{\bf Step 5.} Let $u \in NG$. Then $\wa u \in$ ${\bf aut} \coe$ is an automorphism which restricts to an automorphism $\delta \in {\bf aut} \cog$, and Lemma \ref{lem_cd2} implies that $\delta = \wa y$, where $y := p_*(u) \in QG$. This proves that the diagram (\ref{diag_aut}) commutes.
\end{proof}

As proved in Sec.\ref{intro} for the case $\mB = \coe$, $\mA = \cog$, the vertical arrows of (\ref{diag_aut}) in general are not surjective.

Let $G_1 , G_2$, $S$ be groups endowed with group morphisms $r_1 : G_1 \to S$, $r_2 : G_2 \to S$. Then, it is possible to define the {\em fibered product}
\[
G_1 \times_S G_2
:=
\left\{
( g_1 , g_2 )
\in
G_1 \times G_2 
\ : \
r_1(g_1) = r_2(g_2)
\right\}
\ .
\]

\begin{thm}
\label{thm_stru_autb}
Let $\mA$ be a \sC algebra, $( \rho , \eps , d , S)$ a quasi-special endomorphism, and $\mB$ a mCP of $\mA$ by $\rho$ with a group $G \subseteq$ $C(X^\rho,\sud)$, $G \simeq$ ${\bf aut}_{\mA,\rho}\mB$. Then, the map
\begin{equation}
\label{eq_iso}
r : 
{\bf aut}_{\sigma,\eps} ( \mB ; \mA ) 
\to  
{\bf aut}_{\rho,\eps} \mA \times_{QG} NG
\ \ , \ \
r(\beta) := ( \beta |_\mA , r_\mB (\beta) )
\end{equation}
is a group isomorphism.
\end{thm}

\begin{proof}
The previous theorem implies that $r$ is well-defined. If $r(\beta) = ( id_\mA , 1 )$, then by (\ref{eq_bj}) we conclude that $\beta$ restricts to the identity on $j(\coe)$; since $j(\coe)$ and $\mA$ generate $\mB$, we find $\beta = id_\mB$. Thus, $r$ is injective. Moreover, if $( \alpha , u ) \in$ ${\bf aut}_{\rho,\eps} \mA \times_{QG} NG$, then we have a commutative diagram
\begin{equation}
\label{32}
\xymatrix{
            \mA
	    \ar[r]^{\alpha}
	 &  \mB
	 \\ \cog
	    \ar[u]^{ \mu }
	    \ar[r]^{\subset}
	 &  \coe
	    \ar[u]_{j \circ \wa u}
}
\end{equation}
By the universality property of $\mB =$ $\mA \rtimes_\mu \wa G$ (see \cite[\S 3.1]{Vas2}), we conclude that there is $\beta \in {\bf aut} \mB$ such that
\[
\beta ( a \cdot j(t) ) = \alpha (a) \cdot j \circ \wa u (t)
\ \ , \ \
a \in \mA \ , \ t \in \coe
\ .
\]
In other terms, $r(\beta) =$ $( \alpha , u )$.
\end{proof}

\begin{cor}
\label{cor_ext}
Let $\alpha \in$ ${\bf aut}_{\rho,\eps} \mA$. Then, there is $\beta \in$ ${\bf aut}_{\sigma,\eps} ( \mB ; \mA )$ such that $\beta |_\mA = \alpha$ if and only if 
\[
r_\mA (\alpha) \in p_* (NG) \ .
\]
\end{cor}

\begin{proof}
Let $u \in NG$ such that $p_* (u) =$ $r_\mA (\alpha)$. By the previous theorem, the pair $( \alpha , u )$ defines an element $\beta$ of ${\bf aut}_{\sigma,\eps} ( \mB ; \mA )$ such that $\beta |_\mA = \alpha$.
\end{proof}

\begin{rem}
If $X^\rho$ reduces to a single point,
every element in ${\bf aut}_{\rho,\eps} \mA$ admits an extension in ${\bf aut}_{\sigma,\eps} ( \mB ; \mA )$.
\end{rem}

\begin{rem}
In the general situation described in corollary \ref{cor_ext}, an extension $\beta$ commutes with the action of $G$
if and only if $r_\mB (\beta) \in C_{NG}(G)$, the centralizer of $G$ in $NG$.
\end{rem}

\section{Appendix: Sections of group bundles}

For a complete exposition of the background needed for the present section, we refer to \cite[I.7]{Ste}.
Let $d \in \bN$, and $K \subseteq \sud$ be a closed subgroup with normalizer $NK$ and projection
\begin{equation}
\label{eq_pb}
p_0 : NK \to QK
\ .
\end{equation}
It is well-known that (\ref{eq_pb}) defines a principal bundle with structure group $K$, which in general is nontrivial.
In explicit terms, we may find an open cover $(\Omega_i)$ for $QK$, with homeomorphisms
\[
\pi_i : p_0^{-1} ( \Omega_i )  
\ \to \
\Omega_i \times K
\ .
\]
This implies that the isomorphism class of (\ref{eq_pb}) is uniquely determined by the equivalence class of a cocycle 
\begin{equation}
\label{eq_coc}
\alpha := ( \ (\Omega_i) \ , \ (\alpha_{ij})\  ) 
\ \in \
H^1 ( QK , K )
\ ,
\end{equation}
where the continuous maps $\alpha_{ij} : \Omega_i \cap \Omega_j \to K$ are defined by the compositions $\pi_i \circ \pi_j^{-1}$. Now, it is well-known that the following conditions are equivalent:
(1) The bundle (\ref{eq_pb}) is isomorphic to $QK \times K$; (2) The cocycle (\ref{eq_coc}) is trivial; (3) There exists a continuous section $\gamma : QK \to NK$, $p_0 \circ \gamma =$ $id_{QK}$. 
On the converse, the existence of local sections
\[
\gamma_i : \Omega_i \to NK
\ \ , \ \
p_0 \circ \gamma_i = id_{\Omega_i}
\]
is always ensured, in such a way that $\alpha$ is interpreted as the obstruction to obtain a global section extending the maps $\gamma_i$.

\

Our goal is now to use the previous
construction to provide examples of group bundles such that the map (\ref{def_pnq}) is not surjective. 
Let $X$ be a compact Hausdorff space. We define the group bundles $\mG :=$ $X \times K$, $\mNG :=$ $X \times NK$, $\mQG :=$ $X \times QK$, and consider the associated groups of sections  $G$, $NG$, $QG$, with the projection
\[
p : \mNG \to \mQG \ .
\]
Let us consider a section $s \in QG$; then, $s$ may be regarded as a continuous map
\[
s : X \to QK
\ .
\]
We show that $s$ does not necessarily admit a lift $\wt s : X \to NK$ such that $s = p_* \wt s$, where
\[
p_* : NG \to QG
\ \ , \ \ 
p_* s (u) := p \circ s (x)
\ \ , \ \
x \in X
\ .
\]
For this purpose, we define
\[
X_i := s^{-1}(\Omega_i)
\ ,
\]
and the continuous maps
\[
s_i : X_i \to NK
\ \ , \ \
s_i := \gamma_i \circ s |_{X_i}
\ ,
\]
which by construction satisfy the relations
\[
p_* s_i (x) = s (x) 
\ \ , \ \
x \in X_i
\ .
\]
It can be proved that the following relations are satisfied:
\[
s_i (x) s_j (x)^{-1} = \alpha_{ij} \circ s (x)
\ \ , \ \
x \in X_i \cap X_j
\ .
\]
Thus, by defining
\[
s_* \alpha_{ij} (x) 
\ := \
\alpha_{ij} \circ s (x)
\ \ , \ \
x \in X_i \cap X_j
\ ,
\]
we get a cocycle
\[
s_* \alpha := ( (X_i) \ , \ (s_* \alpha_{ij}) ) 
\ \in \
H^1 (X,K)
\ .
\]
By choosing different local sections $\gamma_i$, the equivalence class of $s_* \alpha$ in $H^1(X,K)$ remains unchanged,
thus $s_* \alpha$ is actually an invariant associated with the section $s$. It turns out that {\em existence of a global section $\wt s : X \to NK$ such that $s = p_* \wt s$ is equivalent to triviality of the cocycle $s_* \alpha$}. In particular, if we pick $X :=$ $QK$ and $s : QK \to QK$ the identity map, then $s_* \alpha = \alpha \in$ $H^1 (QK,K)$, and $s$ admits a lift to $NG$ if and only if $\alpha$ is trivial.

In this way, explicit examples for which the identity map does not admit a lift can be easily provided 
by considering suitable closed subgroups $K \subseteq \sud$.


\end{document}